\documentclass[12pt]{article}

\begin{document}

\null \vspace*{.5in}

\begin{center}
{\Huge \bf Through a Glass Darkly}
\bigskip \\
Steven G. Krantz\footnote{It is a pleasure to thank David H. Bailey, Jonathan Borwein, Robert Burckel, David Collins, 
Marvin Greenberg, Reece Harris, Deborah K. Nelson, and James S. Walker for many useful remarks and 
suggestions about different drafts of this essay.  Certainly their insights have contributed 
a number of significant improvements.}
\end{center}
\vspace*{.25in}

\openup \jot

\section{Prolegomena}
\vspace*{.25in}

\begin{quote}
\large Education is a repetition of civilization in little. 
\medskip \\
\null \hbox{ \ \ } \hfill \rm --- \normalsize Herbert Spencer
\end{quote}
\vspace*{.15in}
	   
Being a mathematician is like being a manic depressive. One
experiences occasional moments of giddy elation, interwoven
with protracted periods of black despair. Yet this is the life
path that we choose for ourselves. And we wonder why nobody
understands us.

The budding mathematician spends an extraordinarily long period
of study and backbreaking hard work in order to attain the
Ph.D. And that is only an entry card into the profession. It
hardly makes one a mathematician.

To be able to call oneself a mathematician, one must have
proved some good theorems and written some good papers
thereon. One must have given a number of talks on his work,
and (ideally) one should have either an academic job or a job
in the research infrastructure. Then, and only
then, can one hold one's head up in the community and call
oneself a peer of the realm. Often one is thirty years old
before this comes about. It is a protracted period of
apprenticeship, and there are many fallen and discouraged and indeed lost
along the way.

The professional mathematician spends his life thinking about
problems that he cannot solve, and learning from his (repeated
and often maddening) mistakes. That he can very occasionally
pull the fat out of the fire and make something worthwhile of
it is in fact a small miracle. And even when he can pull off
such a feat, what are the chances that his peers in the
community will toss their hats in the air and proclaim him a
hail fellow well met? Slim to none at best.

In the end we learn to do mathematics because of its intrinsic
beauty, and its enduring value, and for the personal
satisfaction it gives us. It is an important, worthwhile,
dignified way to spend one's time, and it beats almost any
other avocation that I can think of. But it has its
frustrations.

There are few outside of the mathematical community who have
even the vaguest notion of what we do, or how we spend
our time. Surely they have no sense of what a theorem is, or
how one proves a theorem, or why one would want to.\footnote{From my
solipsistic perspective as a mathematician, this is truly tragic.  For
mathematical thinking is at the very basis of human thought.  It
is the key to an examined life.}  How could
one spend a year or two studying other people's work, only so
that one can spend yet several more years to develop one's own
work?  Were it not for tenure, how could any mathematics ever
get done?

We in the mathematics community expect (as we should) the state
legislature to provide funds for the universities (to pay our
salaries, for instance). We expect the members of Congress to
allocate funds for the National Science Foundation and other
agencies to subvent our research. We expect the White House
Science Advisor to speak well of academics, and of
mathematicians in particular, so that we can live our lives
and enjoy the fruits of our labors. But what do these people
know of our values and our goals? How can we hope that, when
they do the obvious and necessary ranking of priorities that
must be a part of their jobs, we will somehow get sorted near
the top of the list?

This last paragraph explains in part why we as a profession can be
aggravated and demoralized, and why we endure periods of frustration and hopelessness. We
are not by nature articulate---especially at presenting our case to those
who do not speak our language---and we pay a price for that incoherence.
We tend to be solipsistic and focused on our scientific activities, and
trust that the value of our results will speak for themselves. When
competing with the {\tt Wii} and the {\tt iPod}, we are bound therefore to
be daunted. 

\section{Life in the Big City}
\vspace*{.25in}

\begin{quote}
\large The most savage controversies are about those matters 
as to which there is no good evidence either way. \\
\medskip \\
\null \hbox{ \ \ } \hfill \rm --- \normalsize Bertrand Russell
\end{quote}
\vspace*{.15in}

If you have ever been Chair of your department, put in the
position of explaining to the Dean what the department's needs
are, you know how hard it is to explain our mission to the
great unwashed. You waltz into the Dean's office and start
telling him how we must have someone in Ricci flows, we
certainly need a worker in mirror symmetry, and what about
that hot new stuff about the distribution of primes using
additive combinatorics? The Dean, probably a chemist, has no
idea what you are talking about.

Of course the person who had the previous appointment with the
Dean was the Chair of Chemistry, and he glibly told the Dean
how they are woefully shy of people in radiochemistry and
organic chemistry.  And an extra physical chemist or two would be
nice as well.  The Dean said ``sure'', he understood
immediately. It was a real shift of gears then for the Dean to
have to figure out what in the world you (from the Mathematics
Department) are talking about.	How do you put your case in
words that the Dean will understand?  How do you sell yourself
(and your department) to him?\footnote{It is arguable that
a mathematics department is better off with a Dean who
is a musicologist or perhaps a philologist.  Such a scholar
is not hampered by the {\it Realpolitik} of lab science dynamics,
and can perhaps think imaginatively about what our goals are.}

Certainly we have the same problem with society at large.
People understand, just because of their social milieu, why
medicine is important and useful. Computers and their
offspring make good sense; we all encounter computers every
day and have at least a heuristic sense of what they are good
for.  Even certain parts of engineering resonate with the average
citizen (aeronautics, biomedical engineering, civil engineering).
But, after getting out of school, most people have little or no use
for mathematics.  Most financial transactions are handled by machines.
Most of us bring our taxes to professionals for preparation.   Most of
us farm out construction projects around the house to contractors.
If any mathematics, or even arithmetic, is required in the workplace
it is probably handled by software.

One of my wife's uncles, a farmer, once said to me---thinking himself to be in
a puckish mood---that we obviously no longer need mathematicians
because we have computers.  I gave him a patient look and said yes, and
we obviously no longer need farmers because we have vending machines.
He was not amused.  But the analogy is a good one.   Computers are
great for manipulating data, but not for thinking.  Vending machines
are great for handing you a morsel of food {\it that someone else
has produced in the traditional fashion}.

People had a hard time understanding what Picasso's art was about---or even
Andy Warhol's art---but they had a visceral sense that it was interesting
and important. The fact that people would spend millions of dollars for
the paintings gave the activity a certain {\it gravitas}, but there is
something in the nature of art that makes it resonate with our collective
unconscious. With mathematics, people spend their lives coming to grips
with what was likely a negative experience in school, reinforced by
uninspiring teachers and dreadful textbooks. If you are at a cocktail
party and announce that you don't like art, or don't like music, people
are liable to conclude that you are some kind of philistine. If instead you announce
that you don't like mathematics, people conclude that you are a regular
guy. [If you choose to announce that you {\it do} like mathematics, people
are liable to get up and walk away.] To the uninitiated, mathematics is
cold and austere and unforgiving. It is difficult to get even an intuitive
sense of what the typical mathematician is up to. Unlike physicists and
biologists (who have been successfully communicating with the press and
the public for more than fifty years), we are not good at telling
half-truths so that we can paint a picture of our meaning and get our
point across. We are too wedded to the mathematical method. We think in
terms of definitions and axioms and theorems.

\section{Living the Good Life}
\vspace*{.25in}

\begin{quote}
\large One normally thinks that everything that is true is true
for a reason.  I've found mathematical truths that are true for no reason
at all.  These mathematical truths are beyond the power of mathematical 
reasoning because they are accidental and random.  \\
\medskip \\
\null \hbox{ \ \ } \hfill \rm --- \normalsize G. J. Chaitin
\end{quote}
\vspace*{.15in}

The life of a mathematician is a wonderful experience. It is an
exhilarating, blissful existence for those who are prone to enjoy it. One
gets to spend one's time with like-minded people who are in pursuit of a
holy grail that is part of an important and valuable larger picture that
we are all bound to. One gets to travel, and spend time with friends all
over the world, and hang out in hotels, and eat exotic foods, and drink
lovely drinks. One gets to teach bright students and engage in the
marketplace of ideas, and actually to develop new ones. What could be
better? There is hardly a more rewarding way to be professionally engaged.

It is a special privilege to be able to spend one's time---and be
paid for it---thinking original (and occasionally profound) thoughts
and developing new programs and ideas.  One actually feels that
he is changing the fabric of the cosmos, helping people to
see things that they have not seen before, affecting people's lives.\footnote{I have
long been inspired by Freeman Dyson's book [DYS].  It describes both poignantly
and passionately the life of the scientist, and how he can feel that he is
altering and influencing the world around him.}
Teaching can and probably should be a part of this process.  For surely
bringing along the next generation, training a new flank of scholars,
is one of the more enlightened and certainly important pursuits.  Also
interacting with young minds is a beautiful way to stay vibrant and
plugged in, and to keep in touch with the development of new ideas.

Of course there are different types of teaching. The teaching of
rudimentary calculus to freshmen has different rewards from teaching your
latest research ideas to graduate students. But both are important, and
both yield palpable results. What is more, {\it this is an activity that
others understand and appreciate}. If the public does not think of us in
any other way, surely they think of us as teachers. And better that {\it we}
should have to do it.  After all, it is our bailiwick.

The hard fact of the matter is that the powers that be in the university also
appreciate our teaching rather more than they do our many other activities.  After
all, mathematics is a key part of the core curriculum.  A university could
hardly survive without mathematics.  Other majors could not function, could
not advance their students, could not build their curricula, without a basis
in mathematics.  So our teaching role at the institution is both fundamental
and essential.  Our research role is less well understood, especially because
{\it we do not by instinct communicate naturally with scholars in other departments}.

This is actually a key point. We all recall the crisis at the University of
Rochester thirteen years ago, when the Dean shut down the graduate program
in mathematics. His reasoning, quite simply, was that he felt that the
mathematics department was isolated, did not interact productively with
other units on campus, did not carry its own weight. The event at
Rochester rang a knell throughout the profession, for we all knew that
similar allegations could be leveled at any of us. Institutions like
Princeton or Harvard are truly ivory towers, and unlikely to suffer the
sort of indignity being described here. But if you work at a public
institution then look out. I work at a {\it very} private university, and
I can tell you that, in my negotiations as Chair with our Dean, he
sometimes brought up Rochester. And he did {\it not} do so in an effort to
be friendly. He was in fact threatening me.

Some departments, like Earth \& Planetary Science or Biomedical
Engineering, interact very naturally with other subjects. Their
material is intrinsically interdisciplinary. It makes perfect sense
for these people to develop cross-disciplinary curricula and joint majors with
other departments. It is very obvious and sensible for them to apply for
grants with people from departments even outside of their School. A
faculty member of such a department will speak several languages fluently.

It is different for mathematics.  It is a challenge just to speak the one language
of mathematics, and to speak it well.  Most of us do a pretty good job at it, and
those outside of mathematics cannot do it at all.  So there is a natural barrier
to communication and collaboration.   In meetings with other faculty---even from physics
and engineering---we find difficulty identifying a common vocabulary.   We find
that we have widely disparate goals, and very different means of achieving them.

Also our value systems are different.  Our methods for gauging success vary dramatically.
Our reward systems deviate markedly.  Once you become a full Professor you will
serve on tenure and promotion committees for other departments.  This experience is a real
eye-opener, for you will find that the
\vfill
\eject

\noindent criteria used in English and History and Geography are quite
different from what we are accustomed to.\footnote{I still recall serving on the committee
for promotion to Professor of a candidate in Geography.  One of his published writings was
called {\sl A Walk Through China Town}.  It described the experience of walking down Grant
Avenue in San Francisco and smelling the wonton soup.  What would be the analogue
of this in a case for promotion in Mathematics?}  Even our views of truth can be markedly different.

\section{The Why and the Wherefore}
\vspace*{.25in}

\begin{quote}
\large The lofty light of the a priori outshines the dim light of the
world and makes for us incontrovertible truths because of
their ``clearness and distinctness.''
\medskip \\
\null \hbox{ \ \ } \hfill \rm --- \normalsize Ren\'{e} Descartes
\end{quote}
\vspace*{.15in}

A mathematician typically goes through most of his early life as a flaming
success at everything he does. One excels in grade school, one excels in
high school, one excels in college. Even in graduate school one can do
quite well if one is willing to put forth the effort.

Put in slightly different terms:  One can get a long way in the basic material just
by being smart.  Not so much effort or discipline is required.  And this may explain
why so many truly brilliant people get left in the dust.  They reach a point where
some real {\it Sitzfleisch} and true effort are required, and they are simply not
up to it.  They have never had to expend such disciplined study before, so why start now?

While there is no question that being smart can take one a
long way, there comes a point---for {\it all of us}---where it
becomes clear that a capacity for hard work can really make a
difference. Most professional mathematicians put in {\it at least} ten
hours per day, {\it at least} six days per week. There are
many who do much more. And we tend to enjoy it. The great
thing about mathematics is that it does not fight you. It
will not sneak behind your back and bite you. It is always
satisfying and always rewarding.

Doing mathematics is {\it not} like laying bricks or mowing the
grass. The quantity of end product is not a linear function of
the time expended. Far from it. As Charles Fefferman, Fields
Medalist, once said, a good mathematician throws 90\% of his
work in the trash. Of course one learns from all that work,
and it makes one stronger for the next sortie. But one often,
at the end of six months or a year, does not have much to
show.

On the other hand, one can be blessed with extraordinary
periods of productivity. The accumulated skills and insights
of many years of study suddenly begin to pay off, and one
finds that he has plenty to say. And it is {\it quite}
worthwhile. Certainly worth writing up and sharing with others
and publishing. This is what makes life rewarding, and this
is what we live for.

Economists like to use professors as a model, because they run
contrary to many of the truisms of elementary economic theory.
For example, if you pay a Professor of Mathematics twice as
much, that does not mean that he will be able to prove twice
as many theorems, or produce twice as many graduate students.
The truth is that he is probably already working to his
capacity. There are only so many hours in the day. What more
could he do? It is difficult to say what a Professor of
Mathematics should be compensated, because we do not fit the
classical economic model.

Flipped on its head, we could also note that if you give a
Professor of Mathematics twice as much to do, it does not
follow that he will have a nervous breakdown, or quit, or go
into open rebellion. Many of us now have a teaching load of
two courses per semester. But sixty years ago the norm---even
at the very best universities in the United States---was three
courses (or more!) per semester. Also, in those days, there
was very little secretarial help. Professors did a lot of the
drudgery themselves. There were also no NSF grants, and very
little discretionary departmental money, so travel was often
subvented from one's own pocket.  Today life is much better
for everyone.

The fact is that a Professor of Mathematics has a good deal of
slack built into his schedule.  If you double his teaching load,
it means that he has less time to go to seminars, or to talk
to his colleagues, or just to sit and think.  But he will still
get through the day.   Just with considerably less enthusiasm.
And notably less creativity.  Universities are holding faculty
much more accountable for their time these days.  Total Quality
Management is one of many insidious ideas from the business world
that is starting to get a grip at our institutions of higher
learning.  In twenty years we may find that we are much more like
teachers (in the way that we spend our time) and much less like scholars.

Sad to say, the Dean or the Provost has only the vaguest sense of what
our scholarly activities are. When they think of the math department at
all, they think of us as ``those guys who teach calculus.'' They certainly
{\it do not} think of us as ``those guys who proved the Bieberbach
conjecture.'' Such a statement would have little meaning for the typical
university administrator. Of course they are pleased when the faculty
garners kudos and awards, but the awards that Louis de Branges received
for his achievement were fairly low key.\footnote{When I was Chair of the
Mathematics Department, the Dean was constantly reminding me that he
thought of us as a gang of incompetent, fairly uncooperative boobs. One of
his very favorite Chairs at that time was the Head of
Earth \& Planetary Sciences. This man was in fact the leader of the Mars
space probe team, and he actually designed the vehicle that was being used
to explore Mars. Well, you can imagine the kind of presentations that this guy
could give---lots of animated graphics, lots of panoramic vistas, 
lots of dreamy speculation, lots of stories about
other-worldly adventures. His talks were given in the biggest auditoriums
on campus, and they were always packed. The Dean was front and center,
with his tongue hanging out, every time; he fairly glowed in the dark
because he was so pleased and excited. How can a mathematician compete
with that sort of showmanship?  Even if I were to prove the Riemann Hypothesis, it would
pale by comparison.}  They probably would not even raise an eyebrow among the Board
of Trustees.

\section{Such is Life}
\vspace*{.25in}

\begin{quote}
\large There is no religious denomination in which the misuse
of metaphysical expressions has been responsible for so much sin as it has
in mathematics. \\
\medskip \\
\null \hbox{ \ \ } \hfill \rm --- \normalsize Ludwig Wittgenstein
\end{quote}
\vspace*{.15in}

Mathematicians are very much like oboe players.  They do something quite difficult
that nobody else understands.   That is fine, but it comes with a price.

We take it for granted that we work in a rarified stratum of the universe that
nobody else will understand.  We do not expect to be able to communicate with others.
When we meet someone at a cocktail party and say, ``I am a mathematician,'' we expect
to be snubbed, or perhaps greeted with a witty rejoinder like, ``I was never any
good in math.''  Or, ``I was good at math until we got to that stuff with the
letters---like algebra.''

When I meet a brain surgeon I never say, ``I was never any good
at brain surgery. Those lobotomies always got me down.'' When
I meet a proctologist, I am never tempted to say, ``I was
never any good at \dots.'' Why do we mathematicians elicit
such foolish behavior from people?

One friend of mine suggested that what people are really saying to us, when
they make a statement of the sort just indicated, is that they spent
their college years screwing around.  They never buckled down and studied anything
serious.  So now they are apologizing for it.  This is perhaps too simplistic.
For taxi drivers say these foolish things too.  And so do mailmen and butchers.  
Perhaps what people are telling us is that they {\it know} that they should
understand and appreciate mathematics, but they do not.  So instead they
are resentful.

There is a real disconnect when it comes to mathematics.  Most people, by the time
that they get to college, have had enough mathematics so that they can be pretty
sure they do not like it.  They certainly do not want to major in the subject, and their
preference is to avoid it as much as possible.  Unfortunately, for many of these folks,
their major may require a nontrivial amount of math (not so much because the subject
area actually {\it uses} mathematics, but rather because the people who run
the department seem to want to use mathematics as a {\it filter}).  And also unfortunately it happens,
much more often than it should, that people end up changing their majors (from engineering
to psychology or physics to media studies) simply because they cannot hack the math.

In recent years I have been collaborating with plastic surgeons, and I find that
this is a wonderful device for cutting through the sort of conversational impasse
that we have been describing.  {\it Everyone}, at least everyone past a certain
age, is quite interested in plastic surgery.  People want to understand it, they
want to know what it entails, they want to know what are the guarantees of success.
When they learn that there are connections between plastic surgery and mathematics
then that is a hint of a human side of math.  It gives me an entree that I never enjoyed
in the past.

I also once wrote a paper with a picture of the space shuttle in it.  That did not
prove to be quite so salubrious for casual conversation; after all, engineering
piled on top of mathematics does not make the mathematics any more palatable.  But at
least it was an indication that I could speak several tongues.

And that is certainly a point worth pondering if we want to fit into a social
milieu.  Speaking many tongues is a distinct advantage, and gives one a wedge
for making real contact with people.  It provides another way of looking at
things, a new point of contact.  Trying to talk to people {\it about mathematics},
{\it in the language of mathematics}, {\it using the logic of mathematics} is
not going to get you very far.  It will not work with newspaper reporters and it
also will not work with ordinary folks that you are going to meet in the course
of your life.

\section{Mathematics and Art}
\vspace*{.25in}

\begin{quote}
\large  It takes a long time to understand nothing. \\
\medskip \\
\null \hbox{ \ \ } \hfill \rm --- \normalsize Edward Dahlberg 
\end{quote}
\vspace*{.15in}

Even in the times of ancient Greece there was an understanding that mathematics
and art were related.  Both disciplines entail symmetry, order, perspective,
and intricate relationships among the components.  The golden mean is but
one of many artifacts of this putative symbiosis.

M. C. Escher spent a good deal of time at the Moorish castle the Alhambra,
studying the very mathematical artwork displayed there.  This served
to inspire his later studies (which are considered to be a very
remarkable synthesis of mathematics and art).

Today there is more formal recognition of the interrelationship of
mathematics and art. No less an eminence than Louis Vuitton offers a
substantial prize each year for innovative work on the interface of
mathematics and art. Benoit Mandelbrot has received this prize (for his
work on fractals---see [MAN]), and so has David Hoffman for his work with
Jim Hoffman and Bill Meeks on embedded minimal surfaces (see [HOF]).

Mathematics and art make a wonderful and fecund pairing for, as we have
discussed here, mathematics is perceived in general to be austere, unforgiving,
cold, and perhaps even lifeless.  By contrast, art is warm, human, inspiring,
even divine.  If I had to give an after-dinner talk about what I do, I would
not get very far trying to discuss the automorphism groups of pseudoconvex
domains.  I would probably have much better luck discussing the mathematics
in the art of M. C. Escher, or the art that led to the mathematical
work of Celso Costa on minimal surfaces.

Of course we as mathematicians perceive our craft to be an art form.
Those among us who can see---and actually prove!---profound new theorems
are held in the greatest reverence, much as artists.   We see the process
of divining a new result and then determining how to verify it much
like the process of eking out a new artwork.  It would be in our best
interest to convey this view of what we do to the world at large.   Whatever
the merits of fractal geometry may be, Benoit Mandelbrot has done a wonderful
job of conveying both the art and the excitement of mathematics to the
public.  

Those who wish to do so may seek mathematics exhibited in art throughout the ages.  Examples
are
\begin{itemize}
\item  A marble mosaic featuring the small stellated dodecahedron, attributed to Paolo Uccello, in the floor of the San Marco Basilica in Venice.
\item  Leonardo da Vinci's outstanding diagrams of regular polyhedra drawn as illustrations for Luca Pacioli's book {\it The Divine Proportion}.
\item  A glass rhombicuboctahedron in Jacopo de' Barbari's portrait of Pacioli, painted in 1495.
\item  A truncated polyhedron (and various other mathematical objects) which feature in Albrecht D\"{u}rer's engraving Melancholia I.
\item  Salvador Dal\'{\i}'s painting {\sl The Last Supper} in which Christ and his disciples are pictured inside a giant dodecahedron.
\end{itemize}

Sculptor Helaman Ferguson [FER] has made sculptures in various materials of a
wide range of complex surfaces and other topological objects. His work is
motivated specifically by the desire to create visual representations of
mathematical objects.   There are many artists today who conceive of themselves,
and indeed advertise themselves, as mathematical artists.   There are probably
rather fewer mathematicians who conceive of themselves as artistic mathematicians.

Mathematics and music have a longstanding and deeply developed relationship.
Abstract algebra and number theory can be used to understand musical
structure.  There is even a well-defined subject of musical set theory (although
it is used primarily to describe atonal pieces).  Pythagorean tuning is based
on the perfect consonances.  Many mathematicians are musicians, and take great comfort
and joy from musical pastimes.  Music can be an opportunity for mathematicians to interact
meaningfully with a broad cross section of our world.  Mathematicians Noam Elkies and
David Wright have developed wonderful presentations---even full courses---about the symbiosis
between mathematics and music.

Mathematics can learn a lot from art, especially from the way that art reaches
out to humanity.   Part of art is the interface between the artist and the observer.
Mathematics is like that too, but typically the observer is another mathematician.
We would do well, as a profession, to think about how to expand our pool
of observers.

\section{Mathematics vs.\ Physics}
\vspace*{.25in}

\begin{quote}
\large I do still believe that rigor is a relative notion, not an absolute one.
It depends on the background readers have and are expected to use in their judgment. \\
\medskip \\
\null \hbox{ \ \ } \hfill \rm --- \normalsize Ren\'{e} Thom
\end{quote}
\vspace*{.15in}

Certainly ``versus'' is the wrong word here.  Ever since the time
of Isaac Newton, mathematics and physics have been closely allied.
After all, Isaac Newton virtually invented physics as we know it today.
And mathematics in his day was a free-for-all.  So the field was open
for Newton to create any synthesis that he chose.

But mathematics and physics are divided by a common goal, which is to
understand the world around us.  Physicists perceive that ``world'' by observing
and recording and thinking.  Mathematicians perceive that ``world'' by looking within
themselves (but see the next section on Platonism vs.\ Kantianism).

And thus arises a difference in styles. The physicist thinks of himself as
an observer, and is often content to describe what he sees. The
mathematician is {\it never} so content. Even when he ``sees'' with utmost
clarity, the mathematician wants to confirm that vision with a proof. This
fact makes us precise and austere and exacting, but it also sets us apart
and makes us mysterious and difficult to deal with.

I once heard Fields Medalist Charles Fefferman give a lecture (to a mixed 
audience of mathematicians and physicists) about the existence of matter.
In those days Fefferman's goal was to prove the existence of matter
from first principles---in an axiomatic fashion.  I thought that this
was a fascinating quest, and I think that some of the other mathematicians
in the audience agreed with me.  But at some point during the talk a frustrated
physicist raised his hand and shouted, ``Why do you need to do this?  All you
have to do is look out the window to see that matter exists!''

Isn't it wonderful?  Different people have different value systems and different
ways to view the very same scientific facts.  If there is a schism between the
way that mathematicians view themselves and the way that {\it physicists} see us, then
there is little surprise that there is such a schism between our view
of ourselves and the way that non-scientists see us.   Most laymen are content
to accept the world phenomenologically---it is what it is.  Certainly
it is not the average person's job to try to dope out why things are the way they are,
or who made them that way.  This all borders on theology, and that is a distinctly painful
topic.   Better to go have a beer and watch a sporting event on the large-screen TV.
This is {\it not} the view that a mathematician takes.

The world of the mathematician is a world that we have built for ourselves.  And
it makes good sense that we have done so, for we need this infrastructure in
order to pursue the truths that we care about.  But the nature of our
subject also sets us apart from others---even from close allies like
the physicists.   We not only have a divergence of points of view, but also
an impasse in communication.   We often cannot find the words to enunciate
what we are seeing, or what we are thinking.

In fact it has taken more than 2500 years for the modern mathematical mode of
discourse to evolve.  Although the history of proof is rather obscure, we know
that the efforts of Thales and Protagoras and Hippocrates and Theaetetus and Plato and Pythagoras
and Aristotle, culminating in Euclid's magnificent {\it Elements}, have given
us the axiomatic method and the language of proof.  In modern times, the work of 
David Hilbert and Nicolas Bourbaki have helped us to sharpen our focus and nail
down a universal language and methodology for mathematics (see [KRA] for a detailed
history of these matters and for many relevant references).  The idea of mathematical
proof is still changing and evolving, but it is definitely part of who we are and
what we believe.
		
The discussion of Platonism and Kantianism in the next section sheds further
light on these issues.

\section{Plato vs.\ Kant}
\vspace*{.25in}

\begin{quote}
\large It is by logic we prove, it is by intuition that we invent.  \\
\medskip \\
\null \hbox{ \ \ } \hfill \rm --- \normalsize Henri Poincar\'{e}
\end{quote}
\vspace*{.15in}

A debate has been festering in the mathematics profession for a
good time now, and it seems to have heated up in the past few
years (see, for instance [DAV]). And the debate says quite a lot about who we are and
how we endeavor to think of ourselves. It is the question of
whether our subject is Platonic or Kantian.

The Platonic view of the world is that mathematical facts have an independent
existence---very much like classical Platonic ideals---and the research
mathematician {\it discovers} those facts---very much like Amerigo Vespucci discovered
America, or Jonas Salk discovered his polio vaccine.  But it should be clearly
understood that, in the Platonic view, mathematical ideas exist in some higher
realm that is independent of the physical world, and certainly independent of
any particular person.  Also independent of time.  The Platonic view poses the
notion that a theorem can be ``true'' before it is proved.

The Kantian view of the world is that the mathematician creates the subject from
within himself.  The idea of set, the idea of group, the idea
of pseudoconvexity, are all products of the human mind.  They do not exist
out there in nature.  We (the mathematical community) have {\it created} them.

My own view is that both these paradigms are valid, and both play a role
in the life of any mathematician.  On a typical day, the mathematician
goes to his office and sits down and thinks.  He will certainly examine mathematical
ideas that already exist, and can be found in some paper penned by some other mathematician.
But he will also cook things up from whole cloth.  Maybe create a new axiom system, or
define a new concept, or formulate a new hypothesis.  These two activities are by
no means mutually exclusive, and they both contribute to the rich broth that is
mathematics.

Of course the Kantian position raises interesting epistemological questions.  Do we
think of mathematics as being created by each individual?  If that is so, then there
are hundreds if not thousands of distinct individuals creating mathematics from within.
How can they communicate and share their ideas?   Or perhaps the Kantian position
is that mathematics is created by some shared consciousness of the aggregate
humanity of mathematicians.  And then is it up to each individual to ``discover''
what the aggregate consciousness has been creating?  Which is starting
to sound awfully Platonic.  Saunders Mac\,Lane [MAC] argues cogently that mathematical
ideas are elicited or abstracted from the world around is.  This is perhaps
a middle path between the two points of view.

The Platonic view of reality seems to border on theism.  For if mathematical truths
have an independent existence---floating out there in the ether somewhere---then who 
created those truths?  And by what means?  Is it some
\vfill
\eject

\noindent higher power, with
whom we would be well-advised to become better acquainted?  

The Platonic view makes us more like physicists.  It would not make much sense
for a physicist to study his subject by simply making things up.  Or cooking
them up through pure cogitation.   For the physicist is supposed to be
describing the world around him.  A physicist like Stephen Hawking, who
is very creative and filled with imagination, is certainly capable
of cooking up ideas like ``black hole'' and ``supergravity'' and ``wormholes'',
but these are all intended to help explain how the universe works.  They
are not like manufacturing a fairy tale.

There are philosophical consequences for the thoughts expressed in the last
paragraph.  Physicists do not feel honor-bound to prove the claims made
in their research papers.  They frequently use other modes of discourse, ranging
from description to analogy to experiment to calculation.  If we mathematicians
are Platonists, describing a world that is ``already out there'', then why
cannot we use the same discourse that the physicists use?  Why do we need
to be so wedded to proofs?

One can hardly imagine an English Professor trying to decide whether his
discipline is Platonic or Kantian.  Nor would a physicist ever waste his
time on such a quest.  People in those disciplines know where the grist of their
mill lives, and what they are about.  The questions do not really make sense for
them.  We are somewhat alone in this quandary, and it is our job to take possession
of it.  If we can.  

It appears that literary critics and physicists are certainly Platonists.
What else could they be?\footnote{Although a physicist may put a finer
point on it and assert that he has no care for a Platonic realm of ideas.
Rather, he wishes to run experiments and ``ask questions of nature.''} It
is unimaginable that they would cook up their subject from within
themselves. Certainly philosophers can and do engage in this discussion, and
they would also be well-equipped (from a strictly intellectual
perspective) to engage in the Platonic vs.\ Kantian debate. But they have
other concerns. This does not seem to be their primary beat.

The article [MAZ] sheds new and profound light on the questions being considered
here.  This is a discussion that will last a long time, and probably will never
come to any clear resolution.

Once again the Platonic vs.\ Kantian debate illustrates the remove that mathematicians
have from the ordinary current of social discourse.  How can the layman identify
with these questions?  How can the layman even care about them?  If I were a real
estate salesman or a dental technician, what would these questions mean to me?

\section{Seeking the Truth}
\vspace*{.25in}

\begin{quote}
\large In what we really understand, we reason but
little. \\
\medskip \\
\null \hbox{ \ \ } \hfill \rm --- \normalsize William Hazlitt
\end{quote}
\vspace*{.15in}

Mathematicians are good at solving problems.  But we have recognized
for a long time that we have a problem with communicating with laymen, with the public
at large, with the press, and with government agencies.  We have made little
progress in solving this particular problem.   What is the difficulty?

Part of the problem is that we are not well-motivated.  It is not entirely
clear what the rewards would be for solving this problem.  But it is also
not clear what the methodology should be.  Standard mathematical argot
will not turn the trick.  Proceeding from definitions to axioms to theorems will,
in this context, fall on deaf ears.  We must learn a new {\it modus operandi}, and
we must learn how to implement it.

This is not something that anyone is particularly good at, and
we mathematicians have little practice in the matter. We have
all concentrated our lives in learning how to communicate {\it with
each other}. And such activity certainly has its own rewards.
But it tends to make us blind to broader issues. It tends to
make us not listen, and not perceive, and not process the
information that we are given.	Even when useful information
trickles through, we are not sure what to do with it.  It does
not fit into the usual infrastructure of our ideas.  We are
not comfortable processing the data.

This is our own fault.  This is how we have trained ourselves, and it is how
we train our students.  We are not by nature open and outreaching.  We are
rather parochial and closed.  We are more comfortable sticking close to home.
And, to repeat a tired adage, we pay a price for this isolation.

\section{Brave New World}
\vspace*{.25in}

\begin{quote}
\large For most wearers of white coats, philosophy is to science as pornography is to sex: 
it is cheaper, easier, and some people seem, bafflingly, to prefer it. Outside of psychology 
it plays almost no part in the functions of the research machine.
\medskip \\
\null \hbox{ \ \ } \hfill \rm --- \normalsize Steve Jones
\end{quote}
\vspace*{.15in}

For the past 2,000 years, mathematicians have enjoyed a sense of keeping to
themselves, and playing their own tune.\footnote{Although it would be
remiss not to note that Archimedes, Newton, and Gauss were public figures,
and very much a part of society.} It has given us the freedom to think our
own thoughts and to pursue our own truths. By not being answerable to
anyone except ourselves, we have been able to keep our subject pure and
insulated from untoward influences.

But the world has changed around us.  Because of the rise of computers, because
of the infusion of engineering ideas into all aspects of life, because
of the changing nature of research funding, we find ourselves not only
isolated but actually cut off from many of the things that we need in order
to prosper and grow.

So it may be time to re-assess our goals, and our milieu, and indeed our
very {\it lingua franca}, and think about how to fit in more naturally
with the flow of life. Every medical student takes a course on medical
ethics. Perhaps every mathematics graduate student should take a course on
communication.  This would include not only good language skills, but
how to use electronic media, how to talk to people with varying (non-mathematical)
backgrounds, how to seek the right level for a presentation, how to select
a topic, and many of the other details that make for effective verbal and visual
skills.  Doing so would strengthen us as individuals, and it would
strengthen our profession. We would be able to get along more effectively
as members of the university, and also as members of society at large.
Surely the benefits would outweigh the inconvenience and aggravation, and
we would likely learn something from the process. But we must train
ourselves (in some instances {\it re}-train ourselves) to be welcoming to
new points of view, to new perspectives, to new value systems. These
different value systems need not be perceived as inimical to our own.
Rather they are complementary, and we can grow by internalizing them.
	 
Mathematics is one of the oldest avenues of human intellectual endeavor
and discourse.  It has a long and glorious history, and in many ways it
represents the best of what we as a species are capable of doing.  We, the mathematics
profession, are the vessels in which the subject lives.  It is up to us to nurture
it and to ensure that it grows and prospers.  We can no longer do this in isolation.
We must become part of the growing and diversifying process that is human development,
and we must learn to communicate with all parts of our culture.  It is in our best
interest, and it is in everyone else's best interest as well.
\bigskip \bigskip \bigskip \\

\noindent {\Large \sc References}

\begin{enumerate}

\item[{\bf [DAV]}]  E. B. Davies, Let Platonism die, {\it Newsletter of the European
Mathematical Society} 64(1007), 24--25.

\item[{\bf [DYS]}]  F. Dyson, {\it Disturbing the Universe}, Basic Books, New York, 2001.

\item[{\bf [FER]}]  H. Ferguson, Sculpture Gallery, \\
\verb@http://www.helasculpt.com/gallery/index.html@.

\item[{\bf [HOF]}]  D. Hoffman, The computer-aided discovery of new embedded minimal
surfaces, {\it Math.\ Intelligencer} 9(1987), 8--21. 

\item[{\bf [KRA]}]  S. Krantz, {\it The Proof is in the Pudding:  A Look at the Changing Nature
of Mathematical Proof}, Springer Publishing, to appear.

\item[{\bf [MAC]}]  S. Mac\,Lane, Mathematical models:  a sketch for the philosophy of
mathematics, {\it American Mathematical Monthly} 88(1981), 462--472.

\item[{\bf [MAN]}]  B. Mandelbrot, {\it The Fractal Geometry of Nature},
Freeman, New York, 1977.

\item[{\bf [MAZ]}]  B. Mazur, Mathematical Platonism and its opposites, \\
\verb@http://www.math.harvard.edu/~mazur/@.

\end{enumerate}
\vspace*{.42in}

\leftline{Department of Mathematics, Washington University in St. Louis, St. Louis, Missouri 63130}

\leftline{\tt sk@math.wustl.edu}

\end{document}